\documentclass[11pt,letterpaper]{article}

\usepackage[utf8]{inputenc}
\usepackage[margin=1in]{geometry}
\usepackage{amsmath,amsthm,amssymb}
\usepackage{mathtools}
\usepackage{mathrsfs}

\theoremstyle{definition}
\newtheorem{definition}{Definition}[section]
\newtheorem*{definitionstar}{Definition}

\theoremstyle{plain}
\newtheorem{theorem}[definition]{Theorem}

\newtheorem{lemma}[definition]{Lemma}

\newtheorem{proposition}[definition]{Proposition}

\newtheorem{observation}[definition]{Observation}
\newtheorem{claim}[definition]{Claim}

\renewcommand{\le}{\leqslant}
\renewcommand{\ge}{\geqslant}
\renewcommand{\leq}{\leqslant}
\renewcommand{\geq}{\geqslant}

\def \eps {\varepsilon}
\def \es {\emptyset}
\def \tri {\triangle}
\def \sm {\setminus}

\def \F {\mathcal{F}}
\def \G {\mathcal{G}}
\def \mH {\mathcal{H}}
\def \S {\mathcal{S}}
\def \X {\mathcal{X}}
\def \Y {\mathcal{Y}}
\def \U {\mathcal{U}}
\def \V {\mathcal{V}}

\def \ce {\coloneqq}

\title{Maximal 3-wise Intersecting Families with Minimum Size:\\ the Odd Case}
\author{
J\'ozsef Balogh\thanks{Department of Mathematics, University of Illinois at Urbana-Champaign, Urbana, Illinois 61801, USA. E-mail: \texttt{\{jobal, cechen4, haoranl8\}@illinois.edu}.}  \thanks{Research is partially supported by NSF Grant DMS-1764123, NSF RTG grant DMS 1937241, Arnold O. Beckman Research Award (UIUC Campus Research Board RB 22000), and the Langan Scholar Fund (UIUC). %Haoran Luo is partially supported by UIUC Campus Research Board RB 22000.
}
\and
Ce Chen\footnotemark[1]
\and
Haoran Luo\footnotemark[1] \thanks{Research is partially supported by UIUC Campus Research Board RB 22000.}
}
\date{}

\begin{document}
\maketitle

\begin{abstract}
    A family $\F$ on ground set $\{1,2,\ldots, n\}$ is \emph{maximal $k$-wise intersecting} if every collection of $k$ sets in $\F$ has non-empty intersection, and no other set can be added to $\F$ while maintaining this property. Erd\H os and Kleitman asked for the minimum size of a maximal $k$-wise intersecting family. Complementing earlier work of Hendrey, Lund, Tompkins and Tran, who answered this question for $k=3$ and large even $n$, we answer it for $k=3$ and large odd $n$. We show that the unique minimum family is obtained by partitioning the ground set into two sets $A$ and $B$ with almost equal sizes and taking the family consisting of all the proper supersets of $A$ and of $B$.

    A key ingredient of our proof is the stability result by Ellis and Sudakov about the so-called $2$-generator set systems.
\end{abstract}

\section{Introduction}

For a positive integer $n$, denote by $[n]$ the set $\{1,2,\ldots, n\}$. Given a set $S$, we write $2^S$ for the power set of $S$, which is the family of all subsets of $S$, and $S^c$ for the complement set of $S$, which is $[n]\sm S$.
Let $\bar{\F} \ce \{F^c : F \in \F\}$. Then, $\bar{\bar{\F}} = \F$ and $|\bar{\F}| = |\F|$.
Denote the symmetric difference of two sets $S$ and $T$ by $S\Delta T$.

\begin{definitionstar}
A family $\F$ of subsets of $[n]$ is \emph{$k$-wise intersecting} if the intersection of every collection of $k$ distinct sets in $\F$ is non-empty. We call $\F$ \emph{maximal $k$-wise intersecting} if $\F$ is $k$-wise intersecting and no set from $2^{[n]} \sm \F$ can be added to $\F$ keeping the $k$-wise intersecting property.
\end{definitionstar}

%In 1974, Erd\H{o}s and Kleitman~\cite{erdos1974extremal} asked for the minimum possible size of a maximal $k$-wise intersecting family. For $k=3$, Hendrey, Lund, Tompkins and Tran~\cite{hendrey2021maximal} solved the problem when $n$ is a sufficiently large even number. We solve the remaining case for $k=3$ when $n$ is a sufficiently large odd number.

\begin{definitionstar}
\emph{A pair of linked cubes on $[n]$} is a set family $\mH$ of the form $\{A: S \subsetneq A \subseteq[n]\} \cup\left\{B: S^c \subsetneq B \subseteq[n]\right\}$ for some  $S \subseteq[n]$. We call a pair of linked cubes $\mH$ \emph{balanced} if $\left\lfloor\frac{n}{2}\right\rfloor \le |S| \le  \left\lceil\frac{n}{2}\right\rceil$. A balanced pair of linked cubes on $[n]$ contains $2^{\lceil n/2 \rceil} + 2^{\lfloor n/2 \rfloor}-3$ sets. \emph{A balanced pair of cubes on $[n]$} is a set family of the form $\{A: A \subseteq S\}\cup\{B: B\subseteq S^c\}$ where $S\subseteq [n]$ with $\lfloor \frac{n}{2}\rfloor \leq |S|\leq \lceil \frac{n}{2}\rceil$. More generally, a \emph{balanced series of $k$ cubes on $[n]$}, denoted by $\F_{n,k}$, is a set family of the form $\cup_{i=1}^k\{A: A \subseteq S_i\}$ where $S_1,\ldots,S_k$ is a partition of $[n]$ with $\left\lfloor\frac{n}{k}\right\rfloor \le |S_i| \le  \left\lceil\frac{n}{k}\right\rceil$ for each $i\in[k]$.
\end{definitionstar}

 In 1974, Erd\H{o}s and Kleitman~\cite{erdos1974extremal} asked for the minimum possible size of a maximal $k$-wise intersecting family.  Hendrey, Lund, Tompkins and Tran~\cite{hendrey2021maximal} proved that when $n$ is a sufficiently large even number, a maximal $3$-wise intersecting family of minimum size is a balanced pair of linked cubes, thus they determined their minimum size. We will prove an analogue theorem for the case when $n$ is a sufficiently large odd number, thus we fully answer the question of Erd\H os and Kleitman~\cite{erdos1974extremal}  for the case $k=3$.

\begin{theorem}\label{thm::main}
    If $n$ is sufficiently large and $\F$ is a maximal $3$-wise intersecting family  of minimum size on ground set $[n]$, then $\F$ is a balanced pair of linked cubes.
\end{theorem}

%\begin{definitionstar}
%\emph{A balanced series of $k$ cubes on $[n]$}, denoted by $\F_{n,k}$, is a set family of the form $\cup_{i=1}^k\{A: A \subseteq S_i\}$ where $S_1,\ldots,S_k$ is a partition of $[n]$ with $\left\lfloor\frac{n}{k}\right\rfloor \le |S_i| \le  \left\lceil\frac{n}{k}\right\rceil$ for each $i$. We call a balanced series of two cubes $\emph{a balanced pair of cubes}$.
%\end{definitionstar}
In order to complete the proof of the case when $n$ is odd, a key ingredient is the following stability result.

\begin{theorem} \label{mlem}
	Let $n=2\ell+1$ be a sufficiently large odd integer. If $\F$ is a maximal $3$-wise intersecting family on $[n]$ of size $|\F| \le 2^{\lfloor n/2\rfloor}+2^{\lceil n/2 \rceil}-3=3\cdot2^\ell-3$, then there exists a balanced pair of cubes $\F_0 = \{A: A\subseteq S\}\cup \{B: B \subseteq S^c \}$ where $S\subseteq [n]$ with $|S|=\left\lfloor n/2 \right\rfloor=\ell$ such that $|\bar{\F} \Delta \F_0| = o(2^\ell)$.
\end{theorem}

In Section~\ref{sec::ES}, we prove Theorem~\ref{mlem} by modifying  ideas of Ellis and Sudakov~\cite{EllisSudakov2011}. In Section~\ref{sec::pro}, we deduce Theorem~\ref{thm::main} from Theorem~\ref{mlem}, using a method similar to the one used in~\cite{hendrey2021maximal}.  In Section~\ref{sec::largek}, we discuss maximal $k$-wise intersecting families, where   $k\geq 4$.
%It is worth mentioning that when $k\geq 4$ and $n$ is a multiple of $k-1$, although we could not determine the structure of the maximal $k$-wise intersecting family with minimum size, a stability result, similar to Theorem~\ref{mlem} can be obtained by modifying the proof of Ellis and Sudakov~\cite{EllisSudakov2011}.

\section{Proof of Theorem~\ref{mlem}}\label{sec::ES}

\begin{definitionstar}
	A set family $\mH\subseteq 2^{[n]}$ is a \emph{$(1-\eps)$-$k$-generator for $[n]$} if all but at most $\eps 2^n$ subsets of $[n]$ can be expressed as a union of at most $k$ disjoint sets of $\mH$.
\end{definitionstar}

Let $n$ be a sufficiently large odd integer and $\F$ be a maximal $3$-wise intersecting family on $[n]$ of size $|\F| \le 2^{\lfloor n/2\rfloor}+2^{\lceil n/2 \rceil}-3=(\sqrt{2}+ 1/\sqrt{2})2^{n/2}-3$. Let $\F^c \ce 2^{[n]}\sm \F$. Then, there exists some small $\eps=o(1)$ such that $$|\F^c|=2^n-|\F|\geq 2^n-\left(\sqrt{2}+ \frac{1}{\sqrt{2}} \right)\cdot 2^{n/2}+3\geq (1-\eps)2^n.$$
As observed in~\cite{hendrey2021maximal}, for every $A\in \F^c$  there exist $B,C\in\F$ such that $B\cap C\subseteq A^c$ since $\F$ is a maximal $3$-wise intersecting family. Notice that $\F$ is closed upwards, i.e.,~if $F\in \F$ and $F\subseteq F'$, then $F'\in \F$. Thus, we may choose $B\cap C=A^c$ with $B\cup C=[n]$, which is equivalent to $A=B^c\cup C^c$ with $B^c\cap C^c=\emptyset$. Hence, every set in $\F^c$ can be expressed as a union of at most two disjoint sets of $\bar{\F}$, which implies that $\bar{\F}$ is a $(1-\eps)$-$2$-generator for $[n]$. Therefore, Theorem~\ref{mlem} is an immediate corollary of the following theorem.

\begin{theorem}\label{thm::modify}
	Assume $n=2\ell+1$ is a sufficiently large odd integer. For every $\eps'>0$, there exists a small $\delta=\delta(\eps')>0$ depending only on $\eps'$ such that the following holds. If $\G\subseteq 2^{[n]}$ is a $(1-\delta)$-$2$-generator for $[n]$ with $|\G|\leq 2^{\lfloor n/2\rfloor}+2^{\lceil n/2 \rceil}=3\cdot 2^\ell$, then there exists a balanced pair of cubes $\F_0 = \{A: A\subseteq S\}\cup \{B: B \subseteq S^c \}$ where $S\subseteq [n]$ with $|S|=\left\lfloor n/2 \right\rfloor=\ell$ such that $|\G\Delta \F_0|\leq \eps' 2^{\lfloor n/2\rfloor}=\eps'2^\ell$.
\end{theorem}

\noindent\textbf{Proof of Theorem~\ref{mlem}. }Let $\eps'>0$ and $n=2\ell+1$ be a sufficiently large odd integer. Let $\delta=\delta(\eps')>0$ be obtained from Theorem~\ref{thm::modify}. Let $\F$ be a maximal $3$-wise intersecting family on $[n]$ of size $|\F| \le 2^{\lfloor n/2\rfloor}+2^{\lceil n/2 \rceil}-3=3\cdot 2^\ell-3$ and $\G\ce \bar{\F}$. Then, $|\G|=|\F|$ and there exists some small $\eps=o(1)$ such that $\G$ is a $(1-\eps)$-$2$-generator for $[n]$. We may assume that $\eps\leq \delta$ since $n$ is sufficiently large. By Theorem~\ref{thm::modify}, there exists a balanced pair of cubes $\F_0 = \{A: A\subseteq S\}\cup \{B: B \subseteq S^c \}$ where $S\subseteq [n]$ with $|S|=\left\lfloor n/2 \right\rfloor=\ell$ such that $|\bar{\F}\Delta\F_0|=|\G\Delta \F_0|\leq \eps' 2^{\lfloor n/2\rfloor}=\eps'2^\ell$.\qed\\

Before starting the proof of Theorem~\ref{thm::modify}, we need some preparations.

\begin{definitionstar}
For a set family $\mH$, the \emph{disjointness graph} $G_\mH$ is the graph with vertex set $\mH$ and edge set $\{\{F_1,F_2\}\subseteq \mH : F_1\cap F_2 =\es\}$. For two (not necessarily disjoint) set families $\mH_1,\mH_2$, the \emph{disjointness bipartite graph} $G_{\mH_1, \mH_2}$ is the bipartite graph with classes $(\mH_1,\mH_2)$, where there is an edge between $F_1 \in \mH_1$ and $F_2 \in \mH_2$ if and only if $F_1 \cap F_2 = \es$.
\end{definitionstar}

\begin{definitionstar}
	Given set systems $\mH_1, \mH_2$ and a bipartite subgraph $B\subseteq G_{\mH_1}$ with bipartition $(\X,\Y)$, we say $B$ (sometimes say $E(B)$) \emph{generates} $\mH_2$ if every $H\in \mH_2$ can be expressed as a disjoint union of some $X\in\X$ and $Y\in\Y$, i.e., every $H\in \mH_2$ corresponds to an edge of $B$.
\end{definitionstar}

\begin{definitionstar}
	For a set system $\mH$ and $i\in [n]$, let $\mH_i^- \ce \{H\in \mH: i\notin H\}$ be the subfamily of sets not containing $i$ and $\mH_i^+ \ce \{H\sm\{i\}: i\in H\in \mH\}$. Note that $|\mH_i^+|+|\mH_i^-|=|\mH|$.
\end{definitionstar}

\begin{lemma}[Proposition 18 in \cite{EllisSudakov2011}]\label{lem::bipartite}
	Let $c>0$. Then, there exists $b=b(c)>0$ such that for every $\mH\subseteq 2^{[n]}$ with $|\mH|\geq c2^{n/2}$, the disjointness graph $G_{\mH}$ can be made bipartite by deleting at most $\frac{(\log\log n)^b}{\log n}|\mH|^2$ edges.
\end{lemma}

In~\cite{EllisSudakov2011}, Ellis and Sudakov proved a variant of Theorem~\ref{thm::modify} for even $n$, which will be needed for the proof of Theorem~\ref{thm::modify}. We state it below.

\begin{theorem}\label{thm::stability n even}
	For every $k\in\mathbb{N}, \eps, \eps'>0$, there exist $n_0=n_0(k,\eps,\eps')$ and $\eta=\eta(k,\eps,\eps')>0$ such that the following holds. If $n\geq n_0$ is a multiple of $k$ and $\G\subseteq 2^{[n]}$ is a $(1-\eps)$-$k$-generator for $[n]$ with $|\G|\leq (1+\eta) |\F_{n,k}|$, then there exists a balanced series of $k$ cubes $\F_1$ such that $|\G\tri \F_1|\leq \eps' |\F_{n,k}|$.
\end{theorem}

\noindent{\bf Proof of Theorem~\ref{thm::modify}.} Let $n=2\ell+1$ be a sufficiently large odd integer. Suppose $\G \subseteq 2^{[n]}$ is a $(1-\delta)$-$2$-generator for $[n]$ with $|\G|\leq 2^{\lfloor n/2\rfloor}+2^{\lceil n/2 \rceil}=3\cdot 2^\ell$. Then, the number of ways to choose at most two disjoint sets (whose unions are different from each other) from $\G$ is at least $(1-\delta)2^n$, by the definition of $(1-\delta)$-$2$-generators. Hence, $|\G|^2\geq (1-\delta)2^n$, which implies that \[|\G|\geq \sqrt{1-\delta}\cdot 2^{n/2}.\] Moreover, let $G_0 \ce G_{\G}$ be the disjointness graph of $\G$, then \[
1+|\G|+e(G_0)\geq (1-\delta)2^n,\] which implies that \[
e(G_0)\geq (1-\delta)2^{2\ell+1}-3\cdot2^\ell-1.\] We conclude that $G_0$ has edge-density
$$\frac{e(G_0)}{\binom{|\G|}{2}}\geq \frac{(1-\delta)2^{2\ell+1}-3\cdot2^\ell-1}{9\cdot 2^{2\ell-1}}\geq \frac{4-5\delta}{9},$$
where the last inequa\-lity comes from $3\cdot2^\ell+1\leq \delta2^{2\ell-1}$.
Applying Lemma~\ref{lem::bipartite} to $\G$ with $c=\sqrt{1-\delta}$, we get that there exists a constant $b>0$ such that we can delete at most
$$
\frac{(\log\log n)^b}{\log n}|\G|^2\leq \frac{(\log\log n)^b}{\log n}9\cdot 2^{2\ell}=\frac{9(\log\log n)^b}{2\log n}2^{2\ell+1}
$$
edges from $G_0$ and obtain a bipartite graph $G=(\X,\Y)$ with $\X\cup\Y=\G$. Note that $G$ generates all but at most
$$
\delta2^{2\ell+1}+1+|\G|+\frac{9(\log\log n)^b}{2\log n}2^{2\ell+1}
$$
subsets of $[n]$. Since $\frac{9(\log\log n)^b}{2\log n}=o(1)$ and $|\G|\leq 3\cdot 2^\ell=o(2^{2\ell+1})$, we may assume that $G$ generates all but at most $\eps2^n$ subsets of $[n]$ where $\eps=o(1)$ satisfies \begin{equation}\label{eG}
e(G)\geq (1-\eps)2^{2\ell+1}.	
\end{equation}
In particular, $\G$ is a $(1-\eps)$-$2$-generator for $[n]$.

Let $\alpha \ce |\X|/2^\ell$ and $\beta \ce |\Y|/2^\ell$. Since $|\X| + |\Y| = |\G| \le 3 \cdot 2^l$, we have
$$\alpha+\beta\leq 3.$$
Therefore, we have
$$
\alpha\beta\leq 9/4.
$$
Since $\alpha \beta 2^{2\ell}=|\X||\Y|\geq e(G)\geq (1-\eps)2^{2\ell+1}$, we also have $\alpha\beta\geq 2-2\eps$. Therefore,
\begin{equation}\label{eqn::bound alphabeta}
	1-2\eps<\alpha,\beta<2+2\eps.
\end{equation}
Let
\begin{equation*}
\X(1/3)\ce\{i\in [n]: |\X_i^+|\geq |\X|/3\}\quad \text{\ \ \ and\ \ \ \ \ \ }
\Y(1/3)\ce\{i\in [n]: |\Y_i^+|\geq |\Y|/3\}.
\end{equation*}
We characterize the structure of $A\ce\X(1/3)$ and $B\ce\Y(1/3)$ step by step via the following lemma and a series of claims.

\begin{lemma}\label{lem::size}
	%It cannot happen that  each of the equations holds at the same time
	It cannot happen that all the following equations hold at the same time:
	\begin{equation}\label{eqn::size}
		|\X|, |\Y|=(3/2-o(1))2^\ell,\quad
		|\X_n^+|, |\Y_n^+|=(1-o(1))2^{\ell-1} \quad  \textrm{and} \quad
		|\X_n^-|, |\Y_n^-|=(1-o(1))2^\ell.
	\end{equation}
\end{lemma}

\begin{proof}
	Suppose for a contradiction that~\eqref{eqn::size} holds. Since all but at most $\eps2^n=2\eps \cdot 2^{n-1}$ sets in $[n]$ can be expressed as a union of at most two disjoint sets in $\G$, we have $\G_n^-=\X_n^- \cup \Y_n^-$ is a $(1-2\eps)$-$2$-generator for $[n-1]$. By the assumption that~\eqref{eqn::size} holds, $|\G_n^-|=(2-o(1))2^\ell=(2-o(1))2^{(n-1)/2}$. As $|\F_{n-1,2}|=2\cdot 2^{(n-1)/2}-1$, we can apply Theorem~\ref{thm::stability n even} to $\G_n^-$ and conclude that there exists an equipartition $S_1\cup S_2=[n-1]$ such that each of $\X_n^-, \Y_n^-$ contains at least $(1-o(1))2^\ell$ sets in $2^{S_1}, 2^{S_2}$. Define $\U\ce\{F\in \X: F\cap S_2=\emptyset\}, \V\ce\{F\in \Y: F\cap S_1=\emptyset\}$.
    Note that $\U_n^- = \X_n^- \cap 2^{S_1}$ and $\V_n^- = \Y_n^- \cap 2^{S_2}$.
    We have $|\U_n^-|, |\V_n^-|=(1-o(1))2^\ell$, implying
    \begin{equation}\label{eqn::size minus}
    	|\X_n^-\sm \U_n^-|, |\Y_n^-\sm \V_n^-|=o(2^\ell).
    \end{equation}

    \noindent Now we prove that \begin{equation*}
    |\X_n^+\sm \U_n^+|, |\Y_n^+\sm \V_n^+|=o(2^\ell).
    \end{equation*} In fact, for every $X\in\X_n^+\sm \U_n^+$, we have $X\cap S_2\neq \emptyset$ by the definition of $\U$, thus $X\cup\{n\}$ is disjoint from at most $2^{\ell-1}$ subsets of $S_2$. Since $|\Y_n^-\sm 2^{S_2}|=|\Y_n^-\sm \V_n^-|=o(2^\ell)$, the set $X\cup\{n\}$ is disjoint from at most $2^{\ell-1}+o(2^\ell)$ sets in $\Y$. Similarly, for every $Y\in\Y_n^+\sm \V_n^+$, the set $Y\cup\{n\}$ is disjoint from at most $2^{\ell-1}+o(2^\ell)$ sets in $\X$. Let $e_n$ be the number of edges $XY\in E(G)$ such that $n\in X\cup Y$. Since $G$ generates all but at most $\eps2^n$ subsets of $[n]$, at least $2^{2\ell}-\eps2^n=(1-2\eps)2^{2\ell}$ sets containing $n$ correspond to edges of $G$, which implies that \begin{equation*}
    e_n\geq (1-2\eps)2^{2\ell}=\left(1-o(1)\right)2^{2\ell}.
    \end{equation*} Denote $\phi\ce |\U_n^+|/|\X_n^+|$ and $\theta\ce |\V_n^+|/|\Y_n^+|$, then $\phi,\theta\in[0,1]$. Combining with~\eqref{eqn::size}, we have \begin{equation*}
    \begin{aligned}
    e_n &\leq |\U_n^+||\Y_n^-|+|\X_n^+\sm\U_n^+|\left(2^{\ell-1}+o(2^\ell)\right)+|\V_n^+||\X_n^-|+|\Y_n^+\sm\V_n^+|\left(2^{\ell-1}+o(2^\ell)\right)\\
    &=\phi(1-o(1))2^{2\ell-1}+(1-\phi)(1-o(1))2^{2\ell-2}+\theta(1-o(1))2^{2\ell-1}+(1-\theta)(1-o(1))2^{2\ell-2}\\
   &=\left(2+\phi+\theta-o(1)\right)2^{2\ell-2}.
    \end{aligned}
    \end{equation*} Hence, \begin{equation*}
    \left(1-o(1)\right)2^{2\ell}\leq e_n\leq \left(2+\phi+\theta-o(1)\right)2^{2\ell-2},\end{equation*} which implies that $\phi,\theta=1-o(1)$. Therefore, $|\X_n^+\sm \U_n^+|=o(|\X_n^+|)=o(2^{\ell}), |\Y_n^+\sm \V_n^+|=o(|\Y_n^+|)=o(2^{\ell})$ as desired.

    \noindent Now we have \begin{equation*}
    	|\X\sm \U|+ |\Y\sm \V|=|\X_n^-\sm \U_n^-|+|\X_n^+\sm \U_n^+|+|\Y_n^-\sm \V_n^-|+|\Y_n^+\sm \V_n^+|=o(2^\ell). \end{equation*}
    For set $F\in\G$, notice that $F\in (\X\sm \U)\cup (\Y\sm \V)$ if and only if $F\in\X$ satisfies $F\cap S_2\neq \emptyset$ or $F\in\Y$ satisfies $F\cap S_1\neq \emptyset$. For $S\subseteq S_1$, if $S\cup \{n\}\in \X$, then $S\in\X_n^+$ by definition. Hence, there are at least $2^{|S_1|}-|\X_n^+|=2^\ell-(1-o(1))2^{\ell-1}=(1+o(1))2^{\ell-1}$ sets $S\subseteq S_1$ satisfying $S\cup \{n\}\notin \X$. Similarly, there are at least $(1+o(1))2^{\ell-1}$ sets $S'\subseteq S_2$ satisfying $S'\cup \{n\}\notin \Y$. Thus, there are at least $(1+o(1))2^{2\ell-2}$ sets of the form $S\cup S'\cup\{n\}\subseteq [n]$ satisfying $S\subseteq S_1, S\cup \{n\}\notin \X$ and $S'\subseteq S_2, S'\cup \{n\}\notin \Y$. Denote this set system by $\S$. Since $G$ generates all but at most $\eps2^n$ subsets of $[n]$, at least $(1+o(1))2^{2\ell-2}-\eps2^{2\ell+1}=(1-o(1))2^{2\ell-2}$ sets in $\S$ correspond to edges of $G$. If $F\in\S$ can be expressed as a disjoint union of $X\in \X$ and $Y \in \Y$, then either $X\cap S_2\neq\emptyset$ or $Y\cap S_1\neq \emptyset$, implying that either $X$ or $Y$ is in $(\X\sm \U)\cup (\Y\sm \V)$. Hence, the number of choices for $F$ is at most $o(2^\ell)|\G|\leq o(2^\ell)\cdot 3\cdot2^\ell\ll(1-o(1))2^{2\ell-2}$, a contradiction.
    %$\leq o(2^{2\ell})$
\end{proof}

\begin{claim}\label{clm::union}
	$A\cup B=[n]$.
\end{claim}

\begin{proof}
	Suppose for a contradiction that $A\cup B\neq [n]$. We may assume without loss of generality that $n\notin A\cup B$.
	%then $n$ is contained in less than $|\X|/3$ sets of $\X$ and less than $|\Y|/3$ sets of $\Y$.
	Let $x\ce|\X_n^+|/|\X|$ and $y\ce|\Y_n^+|/|\Y|$, then $x,y<1/3$ by the definitions of $A$ and $B$. Recalling that $e_n\geq (1-2\eps)2^{2\ell}$ is the number of disjoint pairs $X\in\X, Y\in\Y$ such that $n\in X\cup Y$, we have \begin{equation}\label{eqn::derive size}
		\begin{aligned}
		(1-2\eps)2^{2\ell}\leq e_n&\leq |\X_n^+||\Y_n^-|+|\Y_n^+||\X_n^-|= |\X_n^+|\left(|\Y|-|\Y_n^+|\right)+|\Y_n^+|\left(|\X|-|\X_n^+|\right)\\
		&=(x+y-2xy)\alpha\beta 2^{2\ell}.
		\end{aligned}
	\end{equation}

    Define the function $f(x,y)\ce x+y-2xy$. On $0\leq x,y\leq 1/3$, the function $f(x,y)$ attains its maximum value $4/9$, when $x=y=1/3$. Combining with~\eqref{eqn::derive size}, we have
    $$1-2\eps\leq f(x,y)\alpha\beta\leq \frac{4}{9}\alpha\beta.$$
    Recalling that $\alpha+\beta\leq 3$, it implies that \begin{equation*}
    	\frac{3}{2}-3\sqrt{\frac{\eps}{2}}\leq \alpha,\beta\leq \frac{3}{2}+3\sqrt{\frac{\eps}{2}}.
    \end{equation*}
    Additionally, $\alpha\beta\leq 9/4$ implies that $1-2\eps\leq f(x,y)\alpha\beta\leq \frac{9}{4}f(x,y)$, so \begin{equation*}
    	\frac{1}{3}-\frac{8\eps}{3}\leq x,y\leq \frac{1}{3}.
    \end{equation*}

    In summary, we have \begin{equation*}
    	|\X|, |\Y|=(3/2-o(1))2^\ell,\quad
    	|\X_n^+|, |\Y_n^+|=(1-o(1))2^{\ell-1} \quad  \textrm{and} \quad
    	|\X_n^-|, |\Y_n^-|=(1-o(1))2^\ell,
    \end{equation*} where $\eps=o(1)$. Therefore, Claim~\ref{clm::union} follows from Lemma~\ref{lem::size}.\end{proof}

\begin{claim}\label{clm::intersection}
	$A\cap B=\emptyset$.
\end{claim}

\begin{proof}
	Suppose for a contradiction that $A\cap B\neq \emptyset$. We may assume without loss of generality that $n\in A\cap B$. Let $x\ce|\X_n^+|/|\X|$ and $y\ce|\Y_n^+|/|\Y|$, then $x,y\geq 1/3$ by the definitions of $A$ and $B$. Notice that \begin{equation*}
		(2-2\eps)2^{2\ell}\leq e(G)\leq |\X||\Y|-|\X_n^+||\Y_n^+|=(1-xy)\alpha\beta 2^{2\ell}.
	\end{equation*}
    Hence, $2-2\eps\leq (1-xy)\alpha\beta\leq \frac{8}{9}\alpha\beta$, which implies that \begin{equation*}
    	\frac{3}{2}-\frac{3}{2}\sqrt{\eps}\leq \alpha,\beta\leq \frac{3}{2}+\frac{3}{2}\sqrt{\eps},
    \end{equation*}
    since $\alpha+\beta\leq 3$. Recalling that $\alpha\beta\leq 9/4$, we have $2-2\eps\leq (1-xy)\alpha\beta\leq \frac{9}{4}(1-xy)$, and hence
    \begin{equation*}
    \frac{1}{3}\leq x,y\leq \frac{1}{3}+\frac{8\eps}{3}.
    \end{equation*}

    In summary, we have \begin{equation*}
    	|\X|, |\Y|=(3/2-o(1))2^\ell,\quad
    	|\X_n^+|, |\Y_n^+|=(1-o(1))2^{\ell-1} \quad  \textrm{and} \quad
    	|\X_n^-|, |\Y_n^-|=(1-o(1))2^\ell,
    \end{equation*} where again $\eps=o(1)$. By Lemma~\ref{lem::size}, we completed the proof of Claim~\ref{clm::intersection}.
\end{proof}

By Claims~\ref{clm::union} and~\ref{clm::intersection}, $A\cup B$ is a partition of $[n]$. We will show that $A\cup B$ is in fact an equipartition of $[n]$ and $\X, \Y$ are not too far from $2^A,2^B$, respectively, thus we prove Theorem~\ref{thm::modify}. The following observation is simple but will be useful hereafter.

\begin{observation}\label{obs}
	If $F\in \X\sm 2^A$, then $F$ has at most $2|\Y|/3$ neighbors in $\Y$. Similarly, if $F\in \Y\sm 2^B$, then $F$ has at most $2|\X|/3$ neighbors in $\X$.
\end{observation}

\begin{proof}
	By symmetry, it suffices to prove the first part. Suppose $F\in \X\sm 2^A$, then there exists $i\in [n]$ such that $i\in F\cap B$. By the definition of $B=\Y(1/3)$, there are at least $|\Y|/3$ sets in $\Y$ containing $i$, which therefore have non-empty intersection with $F$. %Since $G$ is the disjointness graph,
	By the definition of $G$, we conclude that $F$ has at most $|\Y|-|\Y|/3=2|\Y|/3$ neighbors in $\Y$, as desired.
\end{proof}

\begin{claim}\label{clm::size}
	We have $|\X\cap 2^A|\geq (2/3-8\eps)|\X|=(2/3-o(1))|\X|$ and $|\Y\cap 2^B|\geq (2/3-8\eps)|\Y|=(2/3-o(1))|\Y|$. Additionally, $[n]=A\cup B$ is an equipartition.
\end{claim}

\begin{proof}
	Let $\theta \ce \frac{|\X \cap 2^A|}{|\X|}$ and $\phi \ce \frac{|\Y \cap 2^B|}{|\Y|}$. By Observation~\ref{obs} and $\alpha\beta\leq 9/4$, we have\begin{equation*}
			(2-2\eps)2^{2\ell}\leq e(G)\leq |\X\cap 2^A|\cdot|\Y|+|\X\sm 2^A|\cdot\frac{2}{3}|\Y|=\left(\frac{2+\theta}{3}\right)\alpha\beta 2^{2\ell}\leq \frac{9}{4}\left(\frac{2+\theta}{3}\right) 2^{2\ell}.
	\end{equation*}
	Hence, we have $\frac{9}{4}\left(\frac{2+\theta}{3}\right)\geq 2-2\eps$, which implies that $\theta\geq 2/3-8\eps/3=2/3-o(1)$ as desired. Similarly, we can prove $\phi\geq 2/3-8\eps/3= 2/3-o(1)$.

    If $|A|\leq \ell-1$, then
    \begin{equation*}
    	|\X|=\frac{|\X\cap 2^A|}{\theta}\leq \frac{|2^A|}{2/3-8\eps/3}\leq \frac{2^{\ell-1}}{2/3-8\eps/3}<(1-2\eps)2^\ell,
    \end{equation*}
    a contradiction to~\eqref{eqn::bound alphabeta}, so we have $|A|\geq \ell$. Similarly, we have $|B|\geq \ell$. Therefore, $[n]=A\cup B$ is an equipartition.
\end{proof}

Assume without loss of generality that $|A|=\ell$ and $|B|=\ell+1$ from now on. We claim that \begin{equation}\label{eqn::bound2 alphabeta}
	\alpha\leq 3/2+6\eps,\quad \beta\geq 3/2-6\eps,
\end{equation} which are better bounds for $\alpha, \beta$ than~\eqref{eqn::bound alphabeta}. Notice that we only need to show $\beta\geq 3/2-6\eps$, as $\alpha+\beta\leq 3$ will then imply $\alpha\leq 3-\beta\leq 3/2+6\eps$. Suppose that $\beta=3/2-\eta$ with some $\eta\geq 0$, then $\alpha\leq 3/2+\eta$ since $\alpha+\beta\leq 3$. By Observation~\ref{obs}, we have
\begin{equation*}
\begin{aligned}
	&(2-2\eps)2^{2\ell}\leq e(G)\leq |\X\cap 2^A||\Y|+|\X\sm 2^A| \cdot\frac{2}{3}|\Y|= |\X\cap 2^A||\Y|+\left(|\X|-|\X\cap 2^A|\right) \cdot \frac{2}{3}|\Y|\\
	&=\left(\frac{1}{3}|\X\cap 2^A|+\frac{2}{3}|\X|\right)|\Y| \leq \left(\frac{1}{3}\cdot 2^\ell+\frac{2}{3}\left(\frac{3}{2}+\eta\right)2^\ell\right)\left(\frac{3}{2}-\eta\right)2^\ell
	=\left(2-\frac{1}{3}\eta-\frac{2}{3}\eta^2\right)2^{2\ell},
\end{aligned}
\end{equation*} which implies that $2-2\eps\leq 2-\eta/3-2\eta^2/3$. Therefore, $\eta\leq 6\eps$.

The final two claims show that $\X, \Y$ are not too far from $2^A,2^B$, respectively.
\begin{claim}\label{clm::3inES}
	(i) $|2^A\sm \X|\leq 22\eps 2^\ell$.\\
	(ii) $|\Y\sm 2^B|\leq (\sqrt{\eps}+2\eps)2^\ell$.
\end{claim}

\begin{proof}
	Let $$D\ce\min\{|2^B\sm \Y|, |\Y\sm 2^B|\},$$ then $D\leq |\Y\sm 2^B|\leq (1/3+8\eps/3)|\Y|<2|\Y|/3$ by Claim~\ref{clm::size}. Define $\Y'$ from $\Y$ by adding $D$ sets of $2^B\sm \Y$ and deleting $D$ sets of $\Y\sm 2^B$. Thus, $|\Y'|=|\Y|\leq (2+2\eps)2^\ell$ by~\eqref{eqn::bound alphabeta}. Note that $\X$ and $\Y'$ are not necessarily disjoint. Let $G_1 = G_{\X,\Y'}$.
	%Denote   by $G_1$ the disjointness bipartite graph with bipartition $(\X,\Y')$. 
	If $D=|2^B\sm \Y|\leq |\Y\sm 2^B|$, then $2^B\subseteq \Y'$ and $|\Y'\sm 2^B|=|\Y'|-2^{\ell+1}\leq \eps 2^{\ell+1}$; if $D=|\Y\sm 2^B|\leq |2^B\sm \Y|$, then $\Y'\subseteq 2^B$ and $|\Y'\sm 2^B|=0$. In both cases, we have $|\Y'\cap 2^B|\geq|\Y\cap 2^B|$ and \begin{equation}\label{eqn::size difference}
	|\Y'\sm 2^B|\leq \eps 2^{\ell+1}.
	\end{equation} We now compare $e(G_1)$ and $e(G)$. Every deleted $Y\in \Y\sm 2^B$ has at most $2|\X|/3$ neighbors in $\X$ by Observation~\ref{obs}. On the other hand, every added $S\in 2^B\sm\Y$ is disjoint from every set in $2^A\cap \X$, thus has at least $|2^A\cap \X|\geq (2/3-8\eps/3)|\X|$ neighbors in $\X$ by Claim~\ref{clm::size}. Therefore,  \begin{equation*}
	e(G)-e(G_1)\leq D\left(\frac{2}{3}|\X|-\left(\frac{2}{3}-\frac{8\eps}{3}\right)|\X|\right)\leq \frac{2}{3} |\Y|\frac{8\eps}{3}|\X|=\frac{16\eps}{9}\alpha\beta2^{2\ell}\leq \frac{16\eps}{9}\frac{9}{4}2^{2\ell}=2\eps2^{2\ell+1},
	\end{equation*} 
	which, with \eqref{eG}, implies that 
	\begin{equation*}
	e(G_1)\geq e(G)-2\eps2^{2\ell+1}\geq (1-\eps)2^{2\ell+1}-2\eps2^{2\ell+1}=(1-3\eps)2^{2\ell+1}.
	\end{equation*}
	Similarly, we define another bipartite graph $G_2$. Let $$C\ce\min\{|2^A\sm \X|, |\X\sm 2^A|\}.$$ Define $\X'$ from $\X$ by adding $C$ sets of $2^A\sm \X$ and deleting $C$ sets of $\X\sm 2^A$. Thus, $|\X'|=|\X|\geq (1-2\eps)2^\ell$ by~\eqref{eqn::bound alphabeta}. Note that $\X'$ and $\Y'$ are not necessarily disjoint.
	Let $G_2 = G_{\X', \Y'}$.
	%Denote by $G_2$ the disjointness bipartite graph with bipartition $(\X',\Y')$. 
	If $C=|2^A\sm \X|\leq |\X\sm 2^A|$, then $2^A\subseteq \X'$ and $|\X'\cap 2^A|=|2^A|=2^\ell$; if $C=|\X\sm 2^A|\leq |2^A\sm \X|$, then $\X'\subseteq 2^A$ and $|\X'\cap 2^A|=|\X'|\geq (1-2\eps)2^\ell$. In both cases, we have $|\X'\cap 2^A|\geq (1-2\eps)2^\ell$. We now compare $e(G_2)$ and $e(G_1)$. Every deleted $X\in \X\sm 2^A$ intersects $B$, thus has at most $2^\ell$ neighbors in $\Y'\cap 2^B$. By~\eqref{eqn::size difference}, $X$ has at most $\eps2^{\ell+1}$ neighbors in $\Y'\sm 2^B$, thus has at most $(1+2\eps)2^\ell$ in $\Y'$. On the other hand, every added $S\in 2^A\sm\X$ is disjoint from every set in $2^B\cap \Y'$, thus has at least $$|2^B\cap \Y'|=|\Y'|-|\Y'\sm 2^B|\geq (3/2-6\eps)2^\ell-(2\eps)2^\ell=(3/2-8\eps)2^\ell$$ neighbors by~\eqref{eqn::bound2 alphabeta} and~\eqref{eqn::size difference}. Therefore, \begin{equation}\label{eqn::eG2l}
	e(G_2)\geq e(G_1)+C\left(\frac{3}{2}-8\eps-1-2\eps\right)2^\ell\geq (1-3\eps)2^{2\ell+1}+C\left(\frac{1}{2}-10\eps\right)2^\ell.
	\end{equation} If $|\X'|\leq 2^\ell$, then $e(G_2)\leq |\X'||\Y'|\leq |\X'|(3\cdot 2^\ell-|\X'|)$ attains its maximum value $2^{2\ell+1}$ when $|\X'|=2^\ell$. If $|\X'| > 2^l$, then let $|\X'|=(1+\eta)2^\ell$ with some $\eta\geq 0$. By the definition of $\X'$, we have $\X'\supseteq 2^A$. Every $X\in \X'\sm 2^A$ intersects $B$, thus has at most $2^\ell$ neighbors in $\Y'\cap 2^B$, thus has at most $(1+2\eps)2^\ell$ neighbors in $\Y'$ by~\eqref{eqn::size difference}. Since $|\X'|+|\Y'|\leq 3\cdot 2^\ell$, we have $|\Y'|\leq (2-\eta)2^\ell$. Combining with~\eqref{eqn::bound2 alphabeta}, we have $\eta\leq 1/2+6\eps<1$. Therefore, \begin{equation*}
	\begin{aligned}
	e(G_2)&\leq |\X'\sm 2^A|(1+2\eps)2^\ell+|\X'\cap 2^A||\Y'|
	\leq\eta2^\ell(1+2\eps)2^\ell+2^\ell(2-\eta)2^\ell\\
	&=(1+\eps\eta)2^{2\ell+1}\leq (1+\eps)2^{2\ell+1}.
	\end{aligned}
	\end{equation*} In both cases, we have \begin{equation}\label{eqn::eG2u}
	e(G_2)\leq (1+\eps)2^{2\ell+1}. \end{equation} Combining with~\eqref{eG} and~\eqref{eqn::eG2l}, we conclude that \begin{equation}\label{eqn::eG-eG1}
	e(G_1)-e(G)\leq e(G_2)-e(G)\leq (1+\eps)2^{2\ell+1}-(1-\eps)2^{2\ell+1}=\eps2^{2\ell+2}.
\end{equation}
	
	Now we prove (i). By~\eqref{eqn::eG2l} and~\eqref{eqn::eG2u}, we have \begin{equation*}
	C\leq \frac{8\eps}{1/2-10\eps}2^\ell \leq 20\eps 2^\ell,
	\end{equation*}
	so if $C=|2^A\sm \X|$, then we are done. Assume $C=|\X\sm 2^A|$.
	Then, we have
	\begin{equation*}
	|2^A\sm \X|=|2^A|-|2^A\cap \X|=2^{|A|}-(|\X|-|\X\sm 2^A|)\leq 2^\ell-\left((1-2\eps)2^\ell-20\eps 2^\ell\right)=22\eps 2^\ell,
	\end{equation*}
	as desired.
	
	For (ii), we show that it suffices to prove \begin{equation}
	D\leq \sqrt{\eps}2^\ell.
	\end{equation} Indeed, if $D=|\Y\sm 2^B|$, we are done. Assume $D=|2^B\sm \Y|$. Then,
	\begin{equation*}
	|\Y\sm 2^B|=|\Y|-|2^B\cap \Y|=|\Y|-(|2^B|-|2^B\sm \Y|)\leq (2+2\eps)2^\ell-(2-\sqrt{\eps})2^\ell=\left(2\eps+\sqrt{\eps} \right)2^\ell,
	\end{equation*}as desired.
	
	Suppose for a contradiction that $D> \sqrt{\eps}2^\ell$. We claim that there exists some $Y\in \Y\sm 2^B$ having at least $2|\X|/3-8\sqrt{\eps}2^\ell$ neighbors in $\X$. Otherwise, recall that every $S\in \Y'\sm \Y\subseteq 2^B\sm\Y$ has at least $(2/3-8\eps/3)|\X|$ neighbors in $\X$ and $|\X|/3\leq (3/2+6\eps)2^\ell/3=(1/2+2\eps)2^\ell$ by~\eqref{eqn::bound2 alphabeta}. Since $\eps=o(1)$ is small, we have
	\begin{equation*}
	e(G_1)-e(G)\geq D\left(\left(\frac{2}{3}-\frac{8\eps}{3}\right)|\X|-\left(\frac{2|\X|}{3}-8\sqrt{\eps}2^\ell\right)\right)\geq 8\sqrt{\eps}\left(\sqrt{\eps}-\frac{\eps}{2}-2\eps^2\right)2^{2\ell}>\eps 2^{2\ell+2},
	\end{equation*}
	 a contradiction to~\eqref{eqn::eG-eG1}. Suppose some $Y\in \Y\sm 2^B$ has at least $2|\X|/3-8\sqrt{\eps}2^\ell$ neighbors in $\X$. Note that there exists $i\in[n]$ such that $i\in Y\cap A$. We may assume without loss of generality that $i=n$. By the definition of $G$, at most $|\X|-2|\X|/3+8\sqrt{\eps}2^\ell=(1/3+o(1))|\X|$ sets in $\X$ contain $n$, i.e., $|\X_n^+|/|\X| \le 1/3+o(1)$. On the other hand, since $|2^A\sm \X|\leq 22\eps 2^\ell$ by (i), at least $2^{\ell-1}-22\eps 2^\ell=(1-o(1))2^{\ell-1}$ subsets of $A$ containing $n$ are contained in $\X$. Combining with~\eqref{eqn::bound2 alphabeta}, we have $|\X|=(3/2-o(1))2^\ell$, which implies that $|\X_n^+|=(1-o(1))2^{\ell-1}$ and $|\Y|=(3/2-o(1))2^\ell$. By the definition of $B$, we have $|\Y_n^+|/|\Y| \leq 1/3+o(1)$ since $n\in A$. Combining with~\eqref{eqn::derive size}, we get
	 %Moreover, both $\X$ and $\Y$ contain $(1-o(1))2^{\ell-1}$ sets containing $n$, i.e., 
	 $|\Y_n^+|=(1-o(1))2^{\ell-1}$ and hence $|\X_n^-|,|\Y_n^-|=(1-o(1))2^\ell$. Again by Lemma~\ref{lem::size}, we obtain a contradiction and complete the proof of Claim~\ref{clm::3inES}.
\end{proof}

\begin{claim}
	$|2^B\sm \Y|\leq 6\sqrt{\eps}2^\ell$.
\end{claim}

\begin{proof}
	%Denote $2^B\sm \Y$ by $\B$ and $\X\sm 2^A$ by $\X_0$.
	Let $e\in E(G)$. Call $e$ \emph{bad} if $e$ has an endpoint in $\Y\sm 2^B$ and \emph{good} otherwise. By Claim~\ref{clm::3inES} (ii) and~\eqref{eqn::bound2 alphabeta}, $G$ has at most $$|\Y\sm 2^B||\X|\leq \left(\sqrt{\eps}+2\eps\right)2^\ell\cdot \left(3/2+6\eps\right)2^\ell\leq 2\sqrt{\eps}2^{2\ell}$$ bad edges. Fix $S'\in 2^B\sm \Y$ and choose a set $F\subseteq [n]$ of the form $F=S\cup S'$ where $S\subseteq A$.
	%Since $\G$ generates all but at most $\eps2^n$ subsets of $[n]$, at least $2^\ell-\eps 2^n$ subsets of the form $S\cup S'\subseteq [n]$, where $S\subseteq A$, can be expressed as a union of two disjoint sets of $\G$.
	If $F$ corresponds to a good edge of $G$, say $S\cup S'=X\cup Y$ with $X\in\X, Y\in\Y, X\cap Y=\emptyset$, then $Y\in \Y\cap 2^B$ since $XY$ is good. Hence, $Y\subsetneq S'$ is a proper subset of $B$ since $S'\notin\Y$, which implies that $X\in \X\sm 2^A$ since $X$ must intersect $B$. Therefore, we have $X\cap A=S$, which implies that a different $S$ would correspond to a different $X$. There are at most
	$$
	|\X\sm 2^A|=|\X|-|\X\cap 2^A|\leq (3/2+6\eps)2^\ell-(1-22\eps) 2^\ell=(1/2+28\eps)2^\ell
	$$
	different $X$'s by~\eqref{eqn::bound2 alphabeta} and Claim~\ref{clm::3inES} (i). Hence, for fixed $S'\in 2^B\sm \Y$, at most $(1/2+28\eps)2^\ell$ sets $F$ of the form $F=S\cup S'$, where $S\subseteq A$, correspond to good edges of $G$. Recall that $\G$ is a $(1-\eps)$-$2$-generator for $[n]$ and there are $2^\ell|2^B\sm \Y|$ sets of the form $S\cup S'$ with $S\subseteq A, S'\in 2^B\sm \Y$ in total, thus at least $|2^B\sm \Y|(2^\ell-(1/2+28\eps)2^\ell)$ sets of the form $S\cup S'$ with $S\subseteq A, S'\in 2^B\sm \Y$ correspond to bad edges of $G$ or do not correspond to edges of $G$, which implies that \begin{equation*}
    	|2^B\sm \Y|\left(2^\ell-\left(1/2+28\eps\right)2^\ell\right)\leq 2\sqrt{\eps}2^{2\ell}+\eps 2^n=\left(2\sqrt{\eps}+2\eps\right)2^{2\ell}.
    \end{equation*}
    Therefore, we conclude $|2^B\sm \Y|\leq 6\sqrt{\eps}2^\ell$.
\end{proof}

\section{Proof of Theorem~\ref{thm::main}} \label{sec::pro}

Let $n=2\ell+1$ be a sufficiently large odd integer and $\F$ be a maximal $3$-wise intersecting family on $[n]$ of size at most $2^{\lfloor n/2\rfloor}+2^{\lceil n/2 \rceil}-3=3\cdot 2^\ell-3$. Let $\G=\bar{\F}=\{F^c: F\in\F\}$ and fix some $\eps\in(0,1/4)$. By Theorem~\ref{mlem}, there exists $S\subseteq [n]$ of size $\ell$ such that $\F_0\ce\{A: A\subseteq S\}\cup \{B: B\subseteq S^c\}$ satisfies $|\G\Delta\F_0|\leq \eps 2^\ell$.

Recall that $\G$ is closed downwards since $\F$ is closed upwards, so $\emptyset\in\G$. We first prove that $S\notin \G=\bar{\F}$. Suppose for a contradiction that $S^c\in \F$. Since $|\G\Delta\F_0|\leq \eps 2^\ell$, among the $2^{\ell+1}$ subsets of $S^c$, there exists some $A\subseteq S^c$ such that both $A$ and $S^c\sm A$ are contained in $\G$. However, this would imply that $S^c, A^c, (S^c\sm A)^c=S\cup A$ are in $\F$. As $S^c\cap A^c\cap (S\cup A)=\emptyset$, it contradicts that $\F$ is a $3$-wise intersecting family. It can be proved similarly that $S^c\notin \G$.

We work with the following partition $\G=\G_1\cup\G_2\cup\G_3\cup\{\emptyset\}$, where $\G_1\ce\{A\in\G: \emptyset\subsetneq A\subsetneq S\}, \G_2\ce\{B\in\G: \emptyset\subsetneq B\subsetneq S^c\}$ and $\G_3\ce\G\sm\F_0$.

\begin{claim}\label{clm::mpf1}
    $|2^{[n]}\sm (\F\cup \F_0)|\leq |\G_1||\G_2|+|\G_3|\cdot\eps 2^\ell$.
\end{claim}

\begin{proof}
    If $A\notin \F$, then as stated at the beginning of Section~\ref{sec::ES}, there exist $B,C\in \G$ such that $A=B\cup C$ with $B\cap C=\emptyset$. Among those, there exists a pair $\{B,C\}$ such that $$h(\{B,C\})\ce\min \{|B\sm S|+|C\sm S^c|, |B\sm S^c|+|C\sm S|\}$$ attains its minimum value over all such choices of $B$ and $C$. Therefore, we can define an injection $f$ on $2^{[n]}\sm (\F\cup 2^S\cup 2^{S^c})$ by mapping $A$ to a pair of sets $\{B,C\}$ in $\G$, which has minimum $h(\{B,C\})$ and $A=B\cup C, B\cap C=\emptyset$. Since $\emptyset\in \G$, $S,S^c\notin\G$ and $\G$ is closed downwards, we have that $f(A)=\{B,C\}$ must be one of the following two types: $|\{B,C\}\cap \G_1|=|\{B,C\}\cap \G_2|=1$; $|\{B,C\}\cap \G_3|\geq 1$. The number of pairs $\{B,C\}$ satisfying $|\{B,C\}\cap \G_1|=|\{B,C\}\cap \G_2|=1$ is at most $|\G_1||\G_2|$, so it suffices to prove that if $f(A)=\{B,C\}$ with $B\in \G_3$, then there are at most $\eps2^\ell$ choices for $C$.

    Suppose $B=X\cup Y$ where $\emptyset\neq X\subsetneq S, \emptyset\neq Y\subsetneq S^c$. Then, $h(\{B,C\})>0$ and $Y\in \G$, since $Y\subseteq B\in\G$ and $\G$ is closed downwards. There are three possibilities for $C$.

    (1) $C\subseteq S$: Define $B' \ce X\cup C, C'\ce Y$, then $B'\cup C'=B\cup C=A$. Since $h(\{B',C'\})=0<h(\{B,C\})$ and $C'=Y\in\G$, we have that $B'\notin \G$ by the definition of $f$. Note that $B'=X\cup C$ is determined by $C$, so the number of choices for $C\subseteq S$ is at most $|2^S\sm \G|$.

    (2) $C\subseteq S^c$: Similarly, the number of choices for $C$ is at most $|2^{S^c}\sm\G|$.

    (3) $C\in \G\sm\F_0$: The number of choices for such $C$ is at most $|\G\sm\F_0|$.

    In summary, the number of choices for $C$ is at most \begin{equation*}
        |2^S\sm \G|+|2^{S^c}\sm \G|+|\G\sm\F_0|=|\G\Delta\F_0|\leq \eps 2^\ell
    \end{equation*} as desired.
\end{proof}

\begin{claim}\label{clm::mpf2}
    $|\G_1||\G_2|+|\G_3|\cdot\eps 2^\ell\leq (2^\ell-2)(2^{\ell+1}-2)$. Equality holds if and only if $\G_3=\emptyset$.
\end{claim}

\begin{proof}
    By the definitions of $\G_1$ and $\G_2$, we have \begin{equation}\label{eqn::g1g2}
    |\G_1||\G_2|\leq (2^\ell-2)(2^{\ell+1}-2).
    \end{equation} Notice that $\min\{|\G_1|,|\G_2|\}\geq 2^\ell-|\F_0\sm \G|\geq 2^\ell-|\F_0\Delta\G|\geq (1-\eps)2^\ell$. Define the function $g(g_1,g_2,g_3)\ce g_1g_2+g_3\cdot\eps2^\ell$. Given that $\min\{g_1,g_2\}\geq(1-\eps)2^\ell$, we claim that $g(g_1',g_2',g_3')\geq g(g_1,g_2,g_3)$ if $g_1'+g_2'=g_1+g_2+1$ and $g_3'=g_3-1$. In fact, since $g_1'g_2'\in\{(g_1+1)g_2,g_1(g_2+1)\}$, we have \[g(g_1',g_2',g_3')-g(g_1,g_2,g_3)\geq \min\{g_1,g_2\}-\eps2^\ell\geq (1-2\eps)2^\ell>0.\] Therefore, fixing $g_1+g_2+g_3$, the function $g(g_1,g_2,g_3)$ attains its maximum when $g_3=0$. Take $g_i=|\G_i|$ for $i\in[3]$. Combining with~\eqref{eqn::g1g2}, we have \[
    |\G_1||\G_2|+|\G_3|\cdot\eps 2^\ell\leq (2^\ell-2)(2^{\ell+1}-2),\] where equality holds if and only if $\G_3=\emptyset$.
    %By Claim~\refeq{clm::mpf1}, it suffices to prove that as long as there exists some $A\in \G_3$ with $A\notin 2^S\cup 2^{S^c}$, we can replace $A$ by some $B\in \F_0$. In fact, since $|\G|\leq 3\cdot2^\ell-3$, there exists some $B\in \F_0$ such that $B\notin \G$. We replace $A$ by $B$ in $\G$ and consider the change of $|\G_1||\G_2|+|\G_3|\cdot\eps 2^\ell$. Since deleting $A$ from $\G$ would decrease $|\G_1||\G_2|+|\G_3|\cdot\eps 2^\ell$ by at most $\eps 2^\ell$, and adding $B$ to $\G$ would increase $|\G_1||\G_2|+|\G_3|\cdot\eps 2^\ell$ by at least $2^{\ell+1}-|\F_0\sm\G|\geq (1-\eps)2^\ell$, we are done.
\end{proof}

By Claims~\ref{clm::mpf1} and~\ref{clm::mpf2}, we have
\begin{equation*}
    2^{2\ell+1}-(3\cdot2^\ell-1)-|\F|\leq |2^{[n]}\sm (\F\cup \F_0)|\leq (2^\ell-2)(2^{\ell+1}-2)=2^{2\ell+1}-6\cdot2^\ell+4,
\end{equation*}which implies that $|\F|\geq 3\cdot2^\ell-3$. By our assumption that $|\F|\leq 3\cdot2^\ell-3$, equality holds in Claim~\ref{clm::mpf2}. Hence, we have $\G=\F_0$, which means that $\F$ is a balanced pair of cubes.

\section{The Case when $k\geq 4$}\label{sec::largek}

Denote by $f(n,k)$ the minimum possible size of a maximal $k$-wise intersecting family on $[n]$. For the case when $k\geq 4$, Hendrey, Lund, Tompkins and Tran~\cite{hendrey2021maximal} proved the following result.

\begin{proposition}[\cite{hendrey2021maximal}]%prop 4.1
	For $k\geq 4$, there exist positive constants $c_k$ and $d_k$ such that for every positive integer $n$, we have \[
	c_k\cdot2^{n/(k-1)}\leq f(n,k)\leq d_k\cdot2^{n/\lceil k/2\rceil}.\]
\end{proposition}

Very recently, Janzer~\cite{janzer2021note} showed that the lower bound is the right order of magnitude of $f(n,k)$, by constructing a maximal $k$-wise intersecting family of size $O(2^{n/(k-1)})$ for every $k\geq 3$ and $n$. Note that in the special case $k=3$, Theorem~\ref{thm::main} matches Janzer's~\cite{janzer2021note} construction.

Assume that $n$ is sufficiently large. Let $\F$ be a maximal $k$-wise intersecting family on $[n]$ with minimum size. Similarly as in the case $k=3$, one can show that $\bar{\F}$ is a $(1-\eps)$-$(k-1)$-generator for $[n]$, where $\eps=o(1)$. Theorem~\ref{thm::stability n even} and a modification of the method in Section~\ref{sec::pro}
could be used to determine the structure of $\bar{\F}$, if $|\bar{\F} \sm \F_{n,k-1}|$ was small.
%its size was about the same as $|\F_{n,k}|$. 
However, Janzer~\cite{janzer2021note} showed that it is not the case.

\begin{theorem}[\cite{janzer2021note}, Lemma 1.2]\label{janzer1.2}
	%Let $k\geq 4$ and $\F_0$ be a series of $k$ cubes on $[n]$, i.e., there exists a partition $S_1\cup\ldots\cup S_{k-1}$ of $[n]$ with $|S_i|=n/(k-1)+O(1)$ for each $i\in[k-1]$ and $\F_0=2^{S_1}\cup\ldots\cup 2^{S_{k-1}}$. If $|\F\sm \F_0|\leq c\cdot 2^{n/(k-1)}$ and $n$ is large enough compared with $c$, then $\bar{\F}$ cannot be maximal $k$-wise intersecting.
	For every $k \ge 4$, there exist $d= d(k) > 0$, $c = c(k) >0$ and $n_0 = n_0(k) >0$ such that the following holds when $n \ge n_0$. Let $S_1 \cup \ldots \cup S_{k-1}$ be a partition of $[n]$ where $\frac{n}{k-1} -d \le |S_i| \le \frac{n}{k-1} +d$ for every $i \in [k-1]$. Let $\F_0=2^{S_1}\cup\ldots\cup 2^{S_{k-1}}$. For every set family $\F \subseteq 2^{[n]}$ with $|\F \sm \F_0| \le c\cdot 2^{n/(k-1)}$, $\bar{\F}$ cannot be maximal $k$-wise intersecting.
\end{theorem}

Combining our method with Theorems~\ref{thm::stability n even} and~\ref{janzer1.2}, we have the following result.

\begin{proposition}
	For every $k \ge 4$, there exists $c = c(k)>0$ such that  
	$$
	f(n,k)\geq (1+c)|\F_{n,k-1}| = (1+c) \left((k-1)2^{\frac{n}{k-1}} -k + 2\right),
	$$
	when $n$ is divisible by $k-1$.
\end{proposition}
Therefore, a maximum $k$-wise intersecting family is necessary to have a more complex structure when $k \ge 4$.
It is worth mentioning that the exact value of the upper bound on Janzer's construction~\cite{janzer2021note} is $(k-1)2^{k-3}2^{n/(k-1)}-(k-2)(2^{k-1}-1)$, which is larger than $|\F_{n,k-1}|$ with about a multiplicative factor of $2^{k-3}$.

\section*{Acknowledgment}

The authors are grateful for Jingwei Xu, Simon Piga and Andrew Treglown, who participated in fruitful discussions at the  beginning of  the project. Simon and Andrew's visit to University of Illinois was partially supported by NSF RTG grant DMS 1937241.

%\section{Discussions}

%When $n$ is a multiple of $k$, we have the stability but...

%When $n$ is not a multiple of $k$, it is hard to deduce a similar lemma as Lemma~\ref{lem::bipartite}, as ...

\bibliographystyle{abbrv}
\bibliography{reference}

\end{document}